\magnification 1200
          	  	                    \def\R{{\rm I\kern-0.2em R\kern0.2em \kern-0.2em}}
          	  	                    \def\N{{\rm I\kern-0.2em N\kern0.2em \kern-0.2em}}
          	  	                    \def\P{{\rm I\kern-0.2em P\kern0.2em \kern-0.2em}}
          	  	                    \def\B{{\rm I\kern-0.2em B\kern0.2em \kern-0.2em}}
          	  	                    \def\Z{{\rm I\kern-0.2em Z\kern0.2em \kern-0.2em}}
          	  	                    \def\C{{\bf \rm C}\kern-.4em {\vrule height1.4ex width.08em depth-.04ex}\;}
          	  	                    \def\B{{\bf \rm B}\kern-.4em {\vrule height1.4ex width.08em depth-.04ex}\;}
          	  	                    
          	  	                    \def\D{{\Delta}}

          	  	                    \def\z{{\zeta}}

          	  	                    \def\cS{{\cal T}}
          	  	                    \def\cW{{\cal W}}
          	  	                    \def\cV{{\cal V}}

          	  	                    \def\cU{{\cal U}}

          	  	                    \
          	  	                    \vskip 20mm
          	  	                    \centerline {\bf HOLOMORPHIC FUNCTIONS UNBOUNDED ON }
          	  	                    \centerline {\bf CURVES OF FINITE LENGTH}
          	  	                    \vskip 4mm
          	  	                    \centerline{Josip Globevnik}
          	  	                    \vskip 4mm
\noindent\bf Abstract\rm \ \ Given a pseudoconvex domain $D\subset \C^N\ ,N\geq 2$, we prove that 
there is a holomorphic function $f$ on $D$ such that the lengths of paths $p\colon \ [0,1]
\rightarrow D$ along which $\Re f$ is bounded above, with $p(0)$ fixed, grow arbitrarily fast as $p(1)\rightarrow bD$. A 
consequence is the existence of a complete closed complex hypersurface $ M\subset D$ such that
the lengths of paths $p\colon\ [0,1]\rightarrow M$, with $p(0)$ fixed, grow arbitrarily fast as $p(1)\rightarrow bD$. 

\vskip 4mm 
          	  	                    \bf 1.\ Introduction and the main results \rm
          	  	                    \vskip 2mm
          	  	                    Denote by $\D$ the open unit disc in $\C$ and by $\B_N$ the open unit ball 
          	  	                    of $\C^N,\ N\geq 2$. In [G] it was proved that there is a closed complex hypersurface $M$ 
          	  	                    in $\B_N$ which is complete, that is, every path $p\colon [0,1)\rightarrow M$ such that 
          	  	                    $|p(t)|\rightarrow 1$ as $t \rightarrow 1$, has infinite length. This was a consequence of the 
          	  	                    main result of [G], a construction of a holomorphic function on $\B_N$ whose real part 
          	  	                    is unbounded on every path of finite length that ends on $b\B_N$. 
          	  	                    
          	  	                    Recall that a domain $D\subset \C^N,\ N\geq 2$, is \it pseudoconvex\rm \  if it has a 
          	  	                    continuous plurisubharmonic exhaustion function. 
          	  	                    This happens if and only if $D$ is holomorphically convex and if and only if $D$ is 
          	  	                    a domain of holomorphy [H]. Every convex domain is pseudoconvex. 
          	  	                     In the present paper we show that given a pseudoconvex domain $D$ in 
          	  	                    $\C^N,\ N\geq 2$, there is a holomorphic function $f$ on $D$ such that the lengths of  
          	  	                    paths $p\colon [0,1]\rightarrow D$ along which the real part of $f$ is bounded above, 
          	  	                    grow arbitrarily rapidly if $p(0)$ is fixed and $p(1)$ tends to $bD$. 
          	  	                    Our main result is the following
          	  	                    \vskip 2mm
          	  	                    \noindent \bf THEOREM 1.1\ \it Let $D\subset \C^N,\ N\geq 2$, be a pseudoconvex domain 
          	  	                    and let $D_n,\ n\in\N$ be an exhaustion of $D$ by relatively compact open sets
          	  	                    $$
          	  	                    D_1\subset\subset D_2\subset\subset \cdots \subset D,\ \ \ \ \bigcup_{n=1}^\infty D_n = D .
          	  	                    $$
          	  	                    Let $A_n,\ n\in\N$, be an increasing sequence of positive
          	  	                    numbers 
          	  	                    converging to $+\infty$.  There is a function $f$, holomorphic on $D$, with 
          	  	                    the following property:
          	  	                    
          	  	                    \noindent Given $L<\infty$ there is an $n_0\in\N$ such that if $n\geq n_0$
          	  	                    and if $p\colon [0,1]\rightarrow D$ is a path such that

          	  	                    (i)\ \ $\Re [f(p(t))]\leq L\ \ (0\leq t\leq 1)$
          	  	                    
          	  	                    (ii)\  $p(0)\in D_1, \ p(1)\in D\setminus D_n$
          	  	                    
          	  	                    \noindent then the length of $p$ exceeds $A_n$. \rm
          	  	                    \vskip 2mm
          	  	                    \noindent So, in particular, for every  $L<\infty$, the boundary $bD$ 
          	  	                    is infinitely far away for a traveller 
          	  	                    travelling within a sublevel set $\{ z\in D\colon 
          	  	                    \ \Re (f(z))<L\} $ of the real part of $f$ :
          	  	                    \vskip 2mm
          	  	                    \noindent \bf COROLLARY 1.1\ \it 
          	  	                    Given a pseudoconvex domain $D\subset \C^N, \ N\geq 2$, 
          	  	                    there is a holomorphic function on $D$ whose real part is unbounded above on every path 
          	  	                    $p\colon [0,1)\rightarrow D,\ \ p(1)\in bD$, of finite length. \rm
          	  	                    \vskip 2mm
          	  	                    \noindent It is perhaps worth mentioning that for any holomorphic function 
          	  	                    $f$ on $\B_N$ there are paths $p\colon [0,1]\rightarrow \B_N,\ p([0,1))\subset 
          	  	                    \B_N, \ p(1)\in b\B_N$ along which 
          	  	                    $f$ is constant [GS].
          	  	                    
          	  	                    Let $M\subset D$ be a closed complex hypersurface, that is, a closed complex 
          	  	                    submanifold of $D$ of  complex
          	  	                    codimension one. A path $p\colon [0,1)\rightarrow M$ is called \it divergent \rm if 
          	  	                    $p(t)$ leaves every compact subset of $M$ as $t\rightarrow 1$. $M$ is called \it complete \rm if every divergent path $p\colon [0,1)\rightarrow M$
          	  	                    has infinite length. 
          	  	                    
          	  	                    Let $f$ be as in Corollary 1.1. By Sard's theorem one can choose $c\in \C$ such that
          	  	                    $$
          	  	                    M=\{ z\in D\colon\ f(z)=c\} 
          	  	                    $$
          	  	                    is a complex manifold. By the properties of $f$, \ $M$ is a complete hypersurface. So we have
          	  	                    \vskip 2mm
          	  	                    \noindent \bf COROLLARY 1.2\ \it Let $D\subset \C^N,\ N\geq 2,$\ be a pseudoconvex domain. 
          	  	                    Then $D$ contains a complete closed complex hypersurface. \rm
          	  	                    \vskip 2mm
          	  	                    \noindent In the special case when $D=\B_N$, Corollary 1.1 and Corollary 1.2 were proved 
          	  	                    in [G]. In [AL] Corollary 1.2 was proved for 
          	  	                    convex domains in $\C^2$. Note that Theorem 1.1 implies a 
          	  	                     stronger result - given an 
          	  	                    exhaustion $D_j,\ j\in\N$, of a pseudoconvex domain $D$ as in 
          	  	                    Theorem 1.1, there is  
          	  	                    a complete closed complex hypersurface $M$ in $D$ such that along paths, $
          	  	                    M\setminus D_j$ becomes arbitrarily far away as 
          	  	                    $j\rightarrow\infty$: 
          	  	                    
          	  	                   \vskip 2mm
          	  	                    \noindent\bf COROLLARY 1.3\ \it Let $D\subset \C^N,\ N\geq 2,$\ be a pseudoconvex domain,
          	  	                    and let $D_j$ and
          	  	                    $A_j,\ j\in \N$, be as in Theorem 1.1. There is a complete closed complex hypersurface $M$
          	  	                    in $D$ meeting $D_1$ with the following property: There is some $n_0\in\N$ such
          	  	                    that if $n\geq n_0$ then
          	  	                    the length of every 
          	  	                    path $p\colon [0,1]\rightarrow M$, such that $p(0)\in D_1$ and $p(1)\in M\setminus D_n$, 
          	  	                    exceeds $A_n$.
          	  	                
          	  	                    \vskip 4mm
          	  	                    \bf 2.\ The main lemma. Reduction to the case $D=\C^N$\rm 
          	  	                    \vskip 2mm
          	  	                    We assume that $N\geq 2$ 
          	  	                    and write $\B$ for $\B_N$. We shall use \it spherical shells. \rm                  
          	  	If $J$ is an interval contained in $(0,\infty)$ we shall write $\hbox{Sh}(J) =
          	  			    	                    \{ x\in\C^N\colon |x|\in J \}$. So, if $J=[\alpha,\beta ]$ then 
          	  			    	                    \ $\hbox{Sh}(J) =\{ x\in\C^N\colon \alpha\leq |x|\leq \beta\} 
          	  	                    = \beta\overline\B\setminus \alpha\B$. Here is our main lemma:
          	  	                    \vskip 2mm
          	  	\noindent\bf LEMMA 2.1\ \it Let $J=(r,R)$ where $0<r<R<\infty$ and let $A<\infty$.
          	  	                    There is a set $E\subset\hbox{\rm Sh}(J)$ such that
          	  	                    
          	  	                    (i) \ the length of every path $p\colon [0,1]\rightarrow 
          	  	                    \hbox{\rm Sh}(\overline J)\setminus 
          	  	                    E$ such that $|p(0)|=r,\ |p(1)|=R$, exceeds $A$
          	  	                    
          	  	                    (ii) given $\varepsilon >0$ and $L<\infty$ there is a polynomial $\Phi$ on $\C^N$ such 
          	  	                    that $|\Phi|<\varepsilon$ on $r\overline\B$ and $\Re \Phi >L$ on $E$. \rm 
          	  	                    \vskip 2mm
          	  	                    
          	  	                    We will prove Lemma 2.1 in the following sections. To prove Theorem 1.1 we need  
          	  	                    the following consequence of Lemma 2.1. 
          	  	                    \vskip 2mm
          	  	                    \noindent \bf LEMMA 2.2\ \it Let $0<r_1<R_1<r_2<R_2<\cdots ,\ r_n\rightarrow \infty$
          	  			    	                    as $n\rightarrow \infty$, and let $B_n$ be an increasing sequence of positive 
          	  			    	                    numbers converging to $+\infty$. There is a holomorphic function $g$ on $\C^N, 
          	  			    	                    \ N\geq 2$, 
          	  			    	                    such that for any $L<\infty$ there is 
          	  			    	                    an $n_0\in\N$ such that 
          	  			    	                    if $n\geq n_0$ and if $p\colon [0,1]\rightarrow \C^N $ is a path such that 
          	  			    	                    
          	  			    	                    (i)\ $\Re[g(p(t))]\leq L\ \ (0\leq t\leq 1)$
          	  			    	                    
          	  			    	                    (ii) $|p(0)|\leq r_n,\ \ |p(1)|\geq R_n$
          	  			    	                                        
          	  	                    \noindent then the length of $p$ exceeds $B_n$. \rm
          	  	                    \vskip 2mm
          	  \noindent\bf Proof.\rm \ \ Let $0<r_1<R_1<r_2<R_2<\cdots ,\ r_n\rightarrow +\infty $	
          	  as $n\rightarrow \infty$ 
          	  and let $B_n$ be an increasing sequence of positive numbers, converging to $+\infty $.
          	  By Lemma 2.1 there 
          	  is, for each $n$, a set $E_n\subset \hbox{Sh}((r_n,R_n))$ such that
          	  
          	  - the length of every path $p\colon [0,1]\rightarrow \hbox{Sh}([r_n,R_n])\setminus E_n$ 
          	  such that  $|p(0)|=r_n,\ |p(1)|= R_n$, exceeds $B_n$
          	  
          	  - given $\varepsilon >0$ and $L<\infty$ there is a polynomial $\Psi $ on $\C^N$ such 
          	  that $|\Psi|<\varepsilon$ on $r_n\overline\B$ and 
          	  $\Re\Psi >L$ on $E_n$. 
          	  
          	  Let $L_n$ be an increasing sequence converging to $+\infty$. 
          	  Suppose for a moment that we have  a sequence of polynomials $\Phi_n$ such that 
          	  
          	  (a) $\Re \Phi_n > L_n+1$ on $E_n$
          	  
          	  (b) $|\Phi_{n+1}-\Phi_n|<1/2^n$ on $R_n\overline \B$.
          	  
          	  \noindent By (b) the sequence $\Phi_n$ converges uniformly on compacta on $\C^N$, denote by $g$ its limit. 
          	  So $g$ is holomorphic on $\C^N$. On $E_n$ we have  $\Re g = \Re[\Phi_n + \sum_{j=n}^\infty (\Phi_{j+1}-\Phi_j)] \geq \Re \Phi_n - 
          	  \sum_{j=n}^\infty 2^{-j} \geq L_n+1-1 = L_n$. 
          	  
          	  Let $L<\infty$. There is an $n_0$ such that $L<L_n$ for all $n\geq n_0$. Suppose that a path $p$ satisfies (i).
          	  Then there are $\alpha, \ \beta,\ 0<\alpha<\beta <1$ such that $p((\alpha,\beta)
          	  \subset \hbox{Sh}((r_n,R_n))$ and $p(\alpha)=r_n,\ p(\beta)=R_n$. If $p$ satisfies also (i) 
          	  then, since $\Re g\geq L_n $ on $E_n$ it follows that $p|[\alpha,\beta]  $ is a map to 
          	  $\hbox{Sh}([r_n,R_n])\setminus E_n$ so by the preceding discussion the length of
          	  $p|[\alpha _n,\beta_n]$ exceeds $B_n$ and consequently the length of $p$ exceeds $B_n$.
          	  
          	  We construct the sequence 
          	  $\Phi_n$  by induction. Pick a polynomial $\Phi_1$ such that $\Re \Phi_1 > L_1+1$ on $E_1$.
          	  Suppose that we have constructed $\Phi_n$. There 
          	  is a constant $C<\infty$ such that $\Re (\Phi_n+C)\geq L_{n+1}+1$ on $E_{n+1}$. By the preceding discussion there is a polynomial $\Psi$ such that 
          	  $|\Psi |<1/2^n$ on $R_n\overline B$ and $\Re \Psi >C$ . Then $\Phi_{n+1}=\Phi_n+\Psi$ 
          	  satisfies (b) and (a) with $n$ replaced by $n+1$. This completes
          	  the proof.
          	  \vskip 1mm
          	  The same proof gives an analogous result for the ball which 
          	  we will not need in the sequel:
          	  	                    \vskip 2mm
          	  	                    \noindent \bf COROLLARY 2.1\ \it 
          	  	                    Let $0<r_1<R_1<r_2<R_2<\cdots ,\ r_n\rightarrow 1$
          	  			    			    	                    as $n\rightarrow \infty$, 
          	  			 and let $A_n$ be an increasing sequence of positive 
          	  			    			    	                    numbers converging to $\infty$. 
          	  			    			    	                    There is a holomorphic function $g$ on $\B$, 
          	  			    			    	                    such that for any $L<\infty$ there is 
          	  			    			    	                    an $n_0\in\N$ such that 
          	  			    			    	                    if $n\geq n_0$ and if $p\colon [0,1]\rightarrow \overline \B$ 
          	  			    			    	                    is a path such that 
          	  			    			    	                    
          	  			    			    	                    (i)\ $\Re[g(p(t))]\leq L\ \ (0\leq t\leq 1)$
          	  			    			    	                    
          	  			    			    	                    (ii) $|p(0)|\leq r_n,\ \ |p(1)|\geq R_n$
          	  			    			    	                    
          	  			    \noindent then the length of $p$ exceeds $A_n$. \rm
          	  	                    \vskip 2mm 
          	  	                    
          	  	                    Let $D\subset \C^N,\ N\geq 2$, be a 
          	  	                    pseudoconvex domain. Then $D$ is a Stein manifold so there is a 
          	  	                    proper holomorphic embedding  $F\colon D\rightarrow \C^{2N+1}$\ [H, Th.5.3.9].

          	  	                    To prove Theorem 1.1. we first prove the following consequence of 
          	  	                    Lemma 2.2. 
          	  	                    \vskip 2mm
          	  	                    \noindent \bf LEMMA 2.3\ \it Let $0<r_1<R_1<r_2<R_2<\cdots ,\ r_n\rightarrow \infty$
          	  	                    as $n\rightarrow \infty$,   and let $A_n$ be an increasing sequence of positive numbers converging 
          	  	                    to $\infty$. There is a holomorphic function $f$ on $D$ such that for any $L<\infty$ there is 
          	  	                    an $n_0\in\N$ such that 
          	  	                    if $n\geq n_0$ and if $p\colon [0,1]\rightarrow D$ is a path such that

          	  	                    (i) \ $\Re [f(p(t))] \leq L\ \ (0\leq t\leq 1)$
          	  	                    
          	  	                    (ii)  $|F(p(0))|\leq r_n,\ \ |F(p(1))|\geq R_n$
          	  	                    
          	  	                    \noindent then the length of $p$ exceeds $A_n$.\rm 
          	  	                    \vskip 2mm
          	  	                    \noindent Note that Lemma 2.3 implies a more precise version of Theorem 1.1, 
          	  	                    yet for a specific exhaustion $D_n=
          	  	                    \{ z\in D\colon\ |F(z)|<R_n\} ,\ n\in \N$. 
          	  	                    \vskip 2mm \noindent\bf Proof of Lemma 2.3 \rm \ 
          	  	           Let $K_n=
          	  	                    \{ z\in D\colon\ r_n\leq |F(z)|\leq R_n\}.$ 
          	  	                    Let  $p\colon [0,1]\rightarrow K_n$ be a path. Then $q=F\circ p\colon \ 
          	  	                    [0,1]\rightarrow \C^{2N+1}$ is a 
          	  	            path whose length equals
          	  	                    $$
          	  	                    \eqalign{
          	  	                    \hbox{length} (q) &= \int_0^1\biggl|(DF)(p(t))\biggl( {{dp}\over{dt}}(t)\biggr)\biggr| dt\cr 
          	  	                    &\leq \max_{w\in K_n}\|(DF)(w)\| \int_0^1\bigl|{{dp}\over{dt}}(t)\bigr| dt \cr
          	  	                    & = \max_{w\in K_n}\|(DF)(w)\|. \hbox{length}(p) .\cr}
          	  	                    $$
          	  	                    The map $F$ is holomorphic and $K_n$ is compact so
          	  	                    $$
          	  	                    \max_{w\in K_n}\|(DF)(w)\| <\infty .
          	  	                    $$
          	  	                    Let $A_n,\ n\in\N$, be an increasing sequence converging to $+\infty$. Choose an increasing sequence 
          	  	                    $B_n$ converging to $+\infty$ such that 
          	  	                    $$
          	  	                    A_n.\max_{w\in K_n}\|(DF)(w)\| \leq B_n \ (n\in\N) .
          	  	                    \eqno (2.1)
          	  	                    $$
          	  	                    Let $g$ be an entire function on $\C^{2N+1}$ given by Lemma 2.2 and let $f=g\circ F$.
          	  	                    Let $L<\infty$ By Lemma 2.2 there is an $n_0\in\N$ such that if $n\geq n_0$ and if $s\colon 
          	  	                    [0,1]\rightarrow \C^{2N+1}$ is a path such that
          	  	                    
          	  	                    (i) \ $\Re [g(s(t))]\leq L\ \ (0\leq t\leq 1)$
          	  	                    
          	  	                    (ii)  $|s(0)|\leq r_n,\ \ |s(1)| \geq R_n$

          	  	                    \noindent then the length of $s$ exceeds $B_n$.
          	  	                    
          	  	                    Now, let $n\geq n_0$ and let $p\colon [0,1]\rightarrow D$ be a path such that
          	  	                    $|F(p(0))|\leq r_n,\ |F(p(1))|\geq R_n$ and $\Re[f(p(t))]\leq L\ \ (0\leq t\leq 1)$, that is
          	  	                    $$
          	  	                    \Re [g(s(t))]\leq L\ \ (0\leq t\leq 1)
          	  	                    \eqno (2.2)
          	  	                    $$
          	  	                    where $s=F\circ p$. There is a segment $[\alpha, \beta]\subset [0,1]$ such that 
          	  	                    $p|[\alpha,\beta]$ maps $[\alpha,\beta]$ to $K_n$ and $|F(p(\alpha))|=r_n,
          	  	                    \ |F(p(\beta ))|=R_n$, that is, 
          	  	                    $|s(\alpha)|=r_n$,\ $|s(\beta )|=R_n.$  By (2.2) Lemma 2.2 implies that
          	  	                    $$
          	  	                    \hbox{length}(s|[\alpha,\beta]) \geq B_n.
          	  	                    $$
          	  	                    By (2.1) it follows that 
          	  	                    $$
          	  	                    \eqalign{
          	  	                    \hbox{length}(p|[\alpha,\beta]) & \geq { {\hbox{length}((F\circ p)|[\alpha,\beta])}
          	  	                    \over {\max_{w\in K_n}\|(DF)(w)\|}} \cr
          	  	                    & = {{\hbox{length}(s|[\alpha,\beta])}\over{\max_{w\in K_n}\|(DF)(w)\|}} \cr
          	  	                    & \geq {{B_n}\over{\max_{w\in K_n}\|(DF)(w)\|}} \cr
          	  	                    & \geq A_n \cr}
          	  	                    $$
          	  	                    Thus, the length of $p$ exceeds $A_n$. This completes the proof 
          	  	                    of Lemma 2.3 provided that Lemma 2.1 has been proved. 
          	  	                    \vskip 2mm
          	  	                    \noindent \bf Proof of Theorem 1.1. \ \rm Let $D_j,\ j\in\N$, be as in Theorem 1.1 and 
          	  	                    let $w\in D_1$.
          	  	                    Since $F\colon\ D\rightarrow \C^{2N+1}$ is a proper map there are $m_0\in\N$ and 
            a strictly increasing sequence $R_n\nearrow \infty$ such that if $\Delta_n=\{z\in D\colon\ |F(z)|<R_n\} ,\ n\in\N$ 
            then $D_1\subset \Delta_1$ and $\Delta_n\subset D_n\ \ (n\geq m_0)$.
             By Lemma 2.3 there is a holomorphic function $f$ on $D$ such that for any $L<\infty $ there is an $n_0\in\N,\ n_0\geq m_0$, such that 
            if $n\geq n_0$ and if $p\colon [0,1]\rightarrow D$ is a path such that
            $$
             \Re[f(p(t))]\leq L\ \ (0\leq t\leq 1),
            \eqno (2.3)
            $$
            $p(0)\in \Delta_1$ and $p(1)\in D\setminus \D_n$ then the length of $p$ exceeds $A_n$. 
            Since $D_1\subset\Delta_1$ and $\D_n\subset D_n \ (n\geq n_0)$ 
            the same holds for any path $p\colon [0,1]\rightarrow D$ which satisfies (2.3) and 
            $p(0)\in D_1,\  p(1)\in D\setminus D_n$. 
            This completes the proof 
            of Theorem 1.1. 
            \vskip 1mm
            
            \noindent\bf Remark\ \ \rm Note that we used only the fact that $F\colon\ D\rightarrow \C^{2N+1} $
            is a proper holomorphic map. We did not need the fact that it is an injective immersion.
            \vskip 1mm
            It remains to prove Lemma 2.1.
          	  	                    \vskip 4mm
          	  	                    3.\ \bf Proof of Lemma 2.1, Part 1 \rm
          	  	                    \vskip 2mm
          	  	                    If $I_1, I_2$ 
          	  	                    are two intervals contained in $(0,\infty)$ then we shall write 
          	  	                    $I_1<I_2$ provided that $I_1\cap I_2=\emptyset $ and provided that 
          	  	                    $I_2$ is to the right of $I_1$, that is, if $x_1<x_2$ for every $x_1\in I_1,\ x_2\in I_2$. 
          	  	                  \  If $\cV\subset b\B$ is an open set and $J\subset (0,\infty)$ is an interval
          	  	                    then we call 
          	  	       	                  the set 
          	  	       	                  $$K(\cV, J)=\{ tz\colon\ z\in\cV, t\in J\} $$
          	  	       	     a \it spherical box.\rm\ Clearly
          	  	       	     $$
          	  	       	     K(\cV, J)= \{ x\in \hbox{Sh}(J)\colon \ {x\over {|x|}}\in \cV\} .
          	  	       	     $$
          	  	       	     A set of the form $\{ x\in b\B\colon\ 
          	  	                    |x-x_0|<\eta \}$ where $x_0\in b\B$ is called \it a ball in \rm $b\B$ \it of radius \rm
          	  	                    $\eta$.
          	  	                    
          	  	                    \noindent We show that Lemma 2.1 follows from 
          	  	                    \vskip 2mm
          	  	                    \noindent\bf LEMMA 3.1.\ \it There is a $\rho >0$ with the following property : 
          	  	                    
          	  	                    \noindent For every ball $\cV \subset b\B$ of radius $\rho $,\ for every $A<\infty $ and 
          	  	                    for every open interval $J=(\alpha,\beta)$  where $1/2< \alpha<\beta<1$ 
          	  	                    there is a set $E\subset K(\cV, J)$ such that
          	  	                    
          	  	                    (i)\ the length of every path $p\colon [0,1]\rightarrow K(\cV, \overline J)\setminus E$ such that $|p(0)|=\alpha,\ |p(1)|=\beta$, 
          	  	                    exceeds $A$
          	  	                    
          	  	                    (ii) given $\varepsilon >0$ and $L<\infty$ there is a polynomial $\Phi$ on 
          	  	                    $\C^N$ such that $|\Phi|<\varepsilon$ on $\alpha\overline\B $ and $\Re \Phi >L$ on $E$.
          	  	                     \rm
          	  	                    \vskip 2mm
          	  	                    \noindent 
          	  	                    One can view Lemma 3.1 as a local version of Lemma 2.1. Note, however, that $\rho$ 
          	  	                    does not depend on $J$.
          	  	                    \vskip 2mm
          	  	                    \noindent\bf Proof of Lemma 2.1, assuming Lemma 3.1.\  \rm
          	  	                    Let $J=(r,R)$ where $0<r<R<\infty $ and let $A<\infty$.  It is easy to see that it is 
          	  	                    enough to prove Lemma 2.1 in the case when $R$ is close to $r$. Hence, with no loss of generality assume that $1/2<r<R<1$. 
          	  	                    
          	  	                    Let $\rho >0$ be as in Lemma 3.1 and put $\eta =\rho/4$.
          	  	                    Choose points $w_1, w_2,\cdots ,w_M \in b\B$ such that the balls 
          	  	                    $$
          	  	                    \cV_i= \{ w\in b\B\colon \ |w-w_i|<2\eta \} ,\ \ 1\leq i\leq M,
          	  	                    $$ 
          	  	                    cover $b\B$ and then let 
          	  	                  $$
          	  	       	                \cW_i= \{ w\in b\B\colon \ |w-w_i|<4\eta \} ,\ \ 1\leq i\leq M.
          	  	                    $$ 
          	  	                    Then every ball 
          	  	                    $\cV\subset b\B$ of radius $2\eta $ is contained in at least one 
          	  	                    of $\cW_j,\ 1\leq j\leq M$. Thus, if $p\colon [0,1]\rightarrow b\B$ is a path then
          	  	                     either $p([0,1]) $ is contained in $\cW _i$ for some $i,\ 1\leq i\leq M$ 
          	  	                     or else the length of $p$ exceeds $2\eta $. If 
          	  	                     $p\colon [0,1]\rightarrow \hbox{Sh}(\overline J)\subset \hbox{Sh}((1/2,1))$ is a path 
          	  	                     then looking at the radial projections $\pi_R(z)= z/|z|$ we conclude that either $
          	  	         (\pi_R\circ p)([0,1])$ is contained in $\cW_i$ for some $i,\ 1\leq i\leq M$, or 
          	  	         else the length of $\pi_R\circ p$ exceeds $2\eta $ which 
          	  	         implies that the length of $p$ exceeds $\eta$. Thus we have 
          	  	         $$\left. 
          	  	         \eqalign{& 
          	  	         \hbox{Let\ } p\colon [0,1]\rightarrow \hbox{Sh}(J)\ \hbox{be a path.
          	  	         Then either there is some\ } j,\ 1\leq j\leq M \cr
          	  	         &\hbox{such that\ } p([0,1])\subset K(\cW_j,J) \hbox{\ or else the 
          	  	         length of\ } p \hbox{\ exceeds \ }\eta \ .\cr}
          	  	         \right \}
          	  	         \eqno (3.1)
          	  	         $$
          	  	         Choose $\ell \in\N$ so large that 
          	  	         $$
          	  	         \ell \eta >A .
          	  	        \eqno (3.2)
          	  	        $$
          	  	        Divide the interval $J$ into $\ell $ pieces $J_1 = \bigl( r,r+ (R-r)/\ell\bigr)
          	  	        ,\cdots ,
          	  	         J_\nu = (\bigl( r+(\nu-1)(R-r)/\ell, r+\nu (R-r)/\ell \bigr) 
          	  	        , \cdots ,
          	  	        J_\ell = \bigl( r+(\ell -1)(R-r)/\ell , R\bigr)$. 
          	  	        Clearly 
          	  	        $$ (0,r)<J_1<J_2<\cdots <J_\ell < (R,\infty)
          	  	        $$
          	  	        and the length $|J_k|$ of each $J_k,\ 1\leq k\leq \ell$ equals $(R-r)/\ell$. 
          	  	        For each $k,\ 1\leq k\leq \ell$, divide $J_k$ into $M$ equally long pieces 
          	  	        $J_{ks}, 1\leq s\leq M$, so that
          	  	        $J_{ks}$ are pairwise disjoint open intervals contained in $J_k$ such that
          	  	        $$
          	  	        J_{k1}<J_{k2}<\cdots <J_{kM} \hbox{\ \ and\ \ \ } |J_{ks}|= 
          	  	        |J_k|/M = (R-r)/(\ell M) \ \ (1\leq s\leq M).
          	  	        $$
          	  	        For each $k, s,\ 1\leq k\leq \ell, 1\leq s\leq M$, we apply Lemma 3.1 for the ball $\cV =\cW_s$ and the interval $J_{ks}$ 
          	  	        to get a set $E_{ks}\subset K(\cW_s, J_{ks})$ such that 
          	  	        $$\left.
          	  	        \eqalign{
          	  	        &\hbox{if\ } J_{ks}=(\alpha_{ks},\beta_{ks})\ \hbox{then the length of every path\ }
          	  	        p\colon [0,1]\rightarrow K(\cW_s, \overline{J_{ks}})\setminus E_{ks}\cr
          	  	        &\hbox{ such that \ } |p(0)|=\alpha_{ks},\ |p(1)|=\beta_{ks},\hbox{\ exceeds\ }A,\cr }\right\}
          	  	        \eqno (3.3)
          	  	        $$
          	  	        and such that
          	  	        $$
          	  	        \left.\eqalign{&\hbox{given\ }\varepsilon >0  \hbox{\ and\ }L<\infty \hbox{\ there 
          	  	        is a polynomial\ }\Phi\cr
          	  	        &\hbox{on\ }\C^N\hbox{\ such that\ } |\Phi|<\varepsilon \hbox{\ on\ }\alpha_{ks}\overline\B
          	  	        \hbox{\ and\ }\Re\Phi >L\hbox{\ on\ }E_{ks} .\cr}\right\} \eqno (3.4)
          	  	        $$
          	  	        Put
          	  	        $$ E = \bigcup_{s=1}^M\bigcup_{k=1}^\ell E_{ks} .
          	  	        $$
          	  	        We show that $E$ has the required properties. 
          	  	        
          	  	        Clearly $E\subset \hbox{Sh }(J)$. Let $p\colon [0,1]\rightarrow \hbox{Sh}(J)\setminus E$ be a 
          	  	        path such that $|p(0)|=r,\ |p(1)|=R$. Let $J_k = (\alpha_k, \beta_k)\ (1\leq k\leq\ell)$.
          	  	        Since $|p(0)|=r,\ |p(1)|=R$ it follows that for each $k,\ 1\leq k\leq \ell$, there are $\gamma_k, \gamma_k^\prime $
          	  	        such that
          	  	        $$
          	  	        0\leq \gamma_1<\gamma_1^\prime<\gamma_2<\gamma_2^\prime<\cdots<\gamma_\ell <\gamma_\ell^\prime\leq 1, 
          	  	        $$
          	  	        such that for each $k$, $p$ maps $[\gamma_k, \gamma_k^\prime]$ to $\hbox{Sh}(\overline{J_k})$  and 
          	  	        $(\gamma_k, \gamma_k^\prime)$ to $\hbox{Sh}(J_k)$ and satisfies $|p(\gamma_k)| = \alpha_k,\ 
          	  	        |p(\gamma_k^\prime)|=\beta _k$.
          	  	        
          	  	        Fix $k,\ 1\leq k\leq \ell$. By (3.1) there are two possibilities. 
          	  	        Either there is some $s,\ 1\leq s\leq M$, such that $p([\gamma_k, \gamma_k^\prime])\subset K(\cW_s, \overline{J_k})$ 
          	  	        or else the length of $p|[\gamma_k, \gamma_k^\prime]$ exceeds $\eta $. Assume that the first happens. Write 
          	  	        $J_{ks} =(\alpha_{ks},\beta_{ks})$. Since $p([\gamma_k,\gamma_k^\prime])\subset 
          	  	        K(\cW_s,\overline{J_k})$ and 
          	  	        $|p(\gamma_k)|=\alpha_k,\ |p(\gamma_k^\prime )|=\beta_k$ it follows that there are $\gamma_{ks},\ 
          	  	        \gamma_{ks}^\prime$ such that $\gamma_k<\gamma_{ks}<\gamma_{ks}^\prime<\gamma_k^\prime $, 
          	  	        and such that $|p(\gamma_{ks})|=
          	  	        \alpha_{ks},\ p(\gamma_{ks}^\prime)|=\beta_{ks}$ and $p([\gamma_{ks},\gamma_{ks}^\prime ])
          	  	        \subset K(\cW_s, \overline {J_{ks}})$. Clearly $p$ maps $[\gamma_{ks},\gamma_{ks}^\prime ]$ into 
          	  	        $K(\cW_s, \overline{J_{ks}})\setminus E = K(\cW_s, \overline{J_{ks}})\setminus E_{ks}$ 
          	  	        so by (3.3) the length of $p|[\gamma_{ks},\gamma_{ks}^\prime ]$
          	  	        exceeds $A$ and so the length of $p|[\gamma_k,\gamma_k^\prime ]$ exceeds $A$.  This shows that for each $k,\ 1\leq k\leq \ell$,
          	  	        the length of $p|[\gamma_k,\gamma_k^\prime]$ exceeds $\hbox{min}\{ \eta,A\}$ and hence by (3.2) the length of $p$ exceeds $A$. 
          	  	        This shows that $E$ satisfies (i) in Lemma 2.1. 
          	  	        
          	  	        It remains to show (ii) in Lemma 2.1. To this end, rename the intervals $J_{ks},\ 1\leq k\leq\ell,\ 
          	  	        1\leq s\leq M$, into $I_1,I_2,\cdots , I_{\ell M}$ and the sets $E_{ks},\ 1\leq k\leq\ell,\ 
          	  	        1\leq s\leq M$, into $E_1,E_2,\cdots E_{\ell M}$ in such a way that
          	  	        $$
          	  	        (0,r)<I_1<I_2<\cdots <I_{\ell M}<(R,\infty)
          	  	        $$
          	  	        and that 
          	  	        $E_j\subset \hbox{Sh}(I_j)\ \ (1\leq j\leq \ell M).
          	  	        $
          	  	        There are $\mu_j,\ 1\leq j\leq \ell M+1$, such that 
          	  	        $\mu_1=r,\ \mu_{\ell M+1} = R$, and such that 
          	  	        $I_j= (\mu_j, \mu_{j+1})\ \ (1\leq j\leq \ell M)$. Recall that by the properties of $E_j$, 
          	  	        $$
          	  	        \left.
          	  	        \eqalign{
          	  	        &\hbox{for each\ }j,\ 1\leq j\leq \ell M, \hbox{\ and for each\ }\varepsilon>0 \hbox{\ and\ }
          	  	        L<\infty \hbox{\ there\ }\cr
          	  	        &\hbox{is a polynomial\ } \Psi\hbox{\ such that\ }|\Psi|<\varepsilon\hbox{\ on\ } 
          	  	        \mu_j\overline\B\hbox{\ and\ } \Re\Psi > L\hbox{\ on\ }E_j .
          	  	        \cr}
          	  	        \right\}\eqno(3.5)
          	  	        $$
          	  	        Let $L<\infty$ and let $\varepsilon >0$. Let $\Phi_1$ be a polynomial such that
          	  	        $$
          	  	        |\Phi_1|<{\varepsilon\over{\ell M}}\hbox{\ \ on\ \ } \mu_1\overline\B = r\overline B
          	  	        \hbox{\ \ and\ \ }\Re\Phi_1>L+
          	  	        \varepsilon \hbox{\ on\ } E_1
          	  	        $$
          	  	        which is possible by (3.5). We construct polynomials $\Phi_j,\ 2\leq j\leq \ell M$, 
          	  	        such that for each
          	  	        $j,\ 2\leq j\leq \ell M$,
          	  	        $$
          	  	        |\Phi_j-\Phi_{j-1}|<{\varepsilon\over{M\ell }} \hbox{\ \ on\ \ }\mu_j\overline B 
          	  	        \hbox{\ \ and \ \ }
          	  	        \Re \Phi_j> L+\varepsilon \hbox{\ on\ }E_j
          	  	        \eqno (3.6)
          	  	        $$
          	  	        and then put $\Phi = \Phi_{M\ell}$. We show that $\Phi $ has the required 
          	  	        properties. On $r\overline\B =\mu_1\overline B$ we have 
          	  	        $|\Phi|\leq |\Phi_1|+|\Phi_2-\Phi_1|+\cdots +|\Phi_{M\ell}-
          	  	        \Phi_{M\ell -1}|<M\ell.\varepsilon/(M\ell ) =\varepsilon.$ \ \ Fix $j,\ 1\leq j\leq \ell M$. 
          	  	        On $E_j$ we have $\Re \Phi = \Re\bigl[ \Phi_j + (\Phi_{j+1}-\Phi_j)+\cdots +(\Phi_{M\ell}-\Phi_{M\ell -1}\bigr] 
          	  	        \geq \Re \Phi_j - (M\ell -1)\varepsilon/M\ell \geq L+\varepsilon -\varepsilon$. Thus, on 
          	  	        $E=\cup_{j=1}^{M\ell}E_j$ we have $\Re \Phi >L$.
          	  	        
          	  	        To find $\Phi_2,\cdots\Phi_{\ell M} $ satisfying (3.6) we use (3.5): Suppose that we 
          	  	        have constructed $\Phi_j$ where $1\leq j\leq \ell M-1$ . There is a constant 
          	  	        $C<\infty$ such that $\Re \Phi_j+C\geq L+\varepsilon $ on $E_{j+1}$. By (3.5) there 
          	  	        is a polynomial $\Psi $ such that $|\Psi|<\varepsilon/(M\ell )$ on $\mu_{j+1}\overline\B$ and $\Re \Psi>C$ on 
          	  	        $E_{j+1}$. Then $\Phi_{j+1}=\Phi_j+\Psi$ has all the required properties
          	  	        This completes the proof. 
          \vskip 4mm
          \noindent\bf 4.\ Proof of Lemma 3.1, Part 1. \rm
          \vskip 2mm
          Write $M= 2N$ and identify $\C^N$ with $\R^M$ by identifying $(p_1+iq_1,\cdots, p_N+iq_N)\in\C^N$ with 
          $(p_1,q_1,\cdots ,p_N,q_N)\in\R^M$. Let $U_0, U_1, U$ be small open balls in $\R^{M-1}$ centered at 
          the origin, 
          such that 
          $$
          U\subset\subset U_1\subset \subset U_0 .
          $$ 
          Write $W_0=U_0\times (0,\infty)$. This is an open half 
          tube in $\R^M=\C^N$. Similarly,
          write $W_1=U_1\times (0,\infty),\ W=U\times (0,\infty)$. Given an interval $J\subset (1/2,\infty )$ 
          we shall write 
          $$
          W_0(J)=\hbox{Sh}(J)\cap W_0,\ W_1(J)=\hbox{Sh}(J)\cap W_1,\ W(J)=\hbox{Sh}(J)\cap W .
          $$
          
          We assume that the ball $U_0$ is so small that for each $r,\ 1/2<r<1$, the surface $W_0\cap b(r\B)$
          can be written as the graph of the function 
          $$
          \psi_r(x_1,\cdots, x_{M-1}) = \Biggl( r^2-\sum_{i=1}^{M-1}x_i^2\Biggr)^{1/2}
          $$
          defined on $U_0$, that is,
          $$
          W_0\cap b(r\B)=\{ (x_1,\cdots x_{M-1},\psi_r(x_1,\cdots, x_{M-1}))\colon\ (x_1,\cdots ,x_{M-1})\in U_0\} .
          $$
          We now turn to the proof of Lemma 3.1. By rotation it is
          enough to prove that there is one ball
          $\cV\subset b\B$ of
          radius $\rho>0$ with the properties in Lemma 3.1. It is easy to see that to prove this it is enough to 
          prove that there is a ball $U$ as above such that for every $A<\infty$ and for every segment
          $J=(\alpha, \beta]$, \  $1/2<\alpha<\beta<1$, \ there is a set $E\subset W(J)$ such that 
          $$
          \left.
          \eqalign{
          &(i)\hbox{\ the length of every path\ } p\colon\ [0,1]\rightarrow W(\overline J)\setminus E 
          \hbox{\ such that \ }\cr
          &|p(0)|=\alpha,\ |p(1)|= \beta, \hbox{\ exceeds\ } A \cr
          &(ii)\hbox{\ given\ }\varepsilon>0\hbox{\ and\ }L<\infty\hbox{\ there is a polynomial\ }
          \Phi \hbox{\ on\ }\C^N
          \hbox{\ such that\ }\cr
          &|\Phi|<\varepsilon \hbox{\ on\ } \alpha \overline \B 
          \hbox{\ and\ }
          \Re \Phi >L\hbox{\ on \ }E .\cr}\ \ \ \right\} \eqno (4.1)
          $$
          We now use some ideas from [GS] and [G]. Let us describe briefly how the set $E$ will look like.
          We will construct 
          finitely many intervals $J_j,\ 1\leq j\leq n,\ (-\infty,1/2)<J_1<\cdots 
          <J_n < [1,\infty)$. For each $j,\ 
          1\leq j\leq n$, we will construct a convex polyhedral surface $C_j\subset \hbox{Sh}(J_j)$ whose facets will be simplices 
          which is such that $W\setminus C_j$ has two components. From each $C_j$ we 
          shall remove a tiny neighbourhood $\cU_j$ of the 
          skeleton of $C_j$ and what remains intersect with W to get the set $E_j$. The set $E$ will be  the union 
          of $E_j,\ 1\leq j\leq n$. A path $p\colon\ [0,1]\rightarrow W([\alpha, \beta])\setminus E$,\ such that
          $|p(0)|=\alpha,\ |p(1)|=\beta $ will have to pass through each $C_j$, 
          and will have to meet $C_j$ in the neighbourhood $\cU_j$ 
          of $\hbox{Skel}(C_j)$. We shall show that given $A<\infty$, $n\in\N$, the intervals $J_j$, the convex surfaces $C_j$ 
          and $\cU_j,\ 1\leq j\leq n$ can be chosen in such a way that (i) in (4.1) will hold. The fact that $C_j$ are convex 
          and contained in disjoint spherical shells will enable us to satisfy (ii) in (4.1).
          
          Begin with a tessellation $\cS $ of $\R^{M-1}$ into simplices which is 
          periodic with respect to a lattice
          $$\Lambda = \bigl\{ \sum _{i=1}^{M-1} n_ie_i\colon\ n_i\in Z,\ 1\leq i\leq M-1\bigr\}
          $$ where $\{e_1,\cdots ,e_{M-1}\} $ is a basis of $\R^{M-1}$, that is $S+e\in\cS$ 
          for every simplex $S\in\cS$ and for every $e\in\Lambda $. 
          What remains of $\R^{M-1} $ after we remove the interiors of all simplices in $\cS$ 
          we call the \it skeleton\rm\  of $\cS$ and denote by $\hbox{Skel}(\cS)$. More 
          generally, we shall use the tessellations
          $$
          \tau(\cS+z) = \{ \tau(S+z)\colon \ S\in\cS\}
          $$
          where $z\in \R^{M-1}$ and $\tau>0$ and define $\hbox{Skel}(\tau(\cS+z))$ in the same way. 
          
          We now show how to construct the polyhedral surfaces mentioned above. 
          
          Fix $U_0, U_1$ and $U$ as above. Fix $z\in \R^{M-1}$ and let $\tau >0$ be very small. Fix
          $r,\ 1/2<r<1$. To get the vertices of our polyhedral surface we shall "lift" the vertices 
          of each simplex 
          $S\in \tau(\cS+z)$ contained in $U_0$ to $b(r\B)$ in the sense that if $v_1,\cdots, v_M\in\R^{M-1}$ 
          are 
          the vertices of $S$ then $(v_i,\psi_r(v_i)),\ 1\leq i\leq M$ are the vertices of the simplex 
          that we denote by $\Psi_r(S)$. The union of these simplices
          $\Psi_r(S)$ for all 
          $S\in \tau(\cS+z)$ contained 
          in $U_0$ we denote by $\Gamma (r,\tau, z)$. This is a polyhedral surface. 
          It is the graph of 
          the piecewise linear function $\varphi_{r,\tau, z}$ defined on the union
          of all simplices $S$ as 
          above, where, on each such simplex with vertices $v_1,\cdots, v_M$ we have
          $$
          \varphi_{r,\tau, z}\bigl(\sum_{i=1}^M\lambda_i v_i\bigr)= \sum_{i=1}^M \lambda_i\psi_r(v_i)\ \ 
          (0\leq \lambda_i\leq 1,\ 1\leq i\leq M,\  \sum_{i=1}^M\lambda_i=1).
          $$ 
          We will show later that the tessellation $\cS$ can be chosen in such a way that if $U_0$ is chosen small enough 
          then for each $r, \ 1/2<r<1$, and each $z$, the surface $\Gamma(r,\tau, z)$ 
          will be \it convex\rm\ in 
          the sense that given a simplex $\Psi_r(S)$ where $S\in \tau(\cS+z)$ is contained in
          $U_0$, the 
          intersection of the hyperplane $H$ containing $\Psi_r(S)$ with $\Gamma(r,\tau, z)$ is 
          precisely $\Psi_r(S)$, that is, all of $\Gamma (r,\tau, z)$ except $\Psi_r(S)$ is contained in the open 
          halfspace bounded by $H$ which contains the origin.
          
          Let $d$ be the length of the longest edge of the simplices in $\cS$. 
          Then $\tau d$ is the length of the 
          longest edge of the simplices in $\tau(\cS+z)$ for any $\tau>0$ and 
          any $z\in\R^{N-1}$. There is a constant
          $\nu >0$ depending on $U_0$ such that for each $r,\ 1/2<r<1$, the length 
          of the longest 
          edge of a simplex building
          $\Gamma (r,\tau, z)$ does not exceed $\gamma=\nu\tau d$.  Thus the vertices of each such simplex are contained 
          in a spherical cap $\{ x\in b(r\B)\colon\ |x-x_0|<\gamma\}$ for some $x_0\in b(r\B)$ and
          consequently the simplex is contained in the convex hull of this spherical cap. If
          $1/2<r<1$ then it is easy to see that this convex hull misses $(r-2\gamma^2)\overline\B$. It follows that 
           there is a constant
          $\omega =
          2\nu^2d^2 $ such that 
          $$
          \left.\eqalign{
          \hbox{for each\ }r, \ 1/2<r<1, \hbox{\ each\ }z\in\R^{M-1}\hbox{\ and each\ }\cr \tau >0
          \hbox{\ we have\ }\Gamma(r,\tau, z) \subset \hbox{Sh}((r-\omega\tau^2, r]) .\cr}
          \ \ \ \right\}
          \eqno (4.2)
          $$
          It is a simple geometric fact that there is a $\delta>0$ such that
          $$
          \left.\eqalign{
          \hbox{\ for every\ }r,\  1/2<r<1,\hbox{\ and for every \ } s,\ r-\delta <s<r, 
          \hbox{\ every line\ }\cr
          \hbox{\ that meets\ }W((s,r])\hbox{\ and misses\ }W_1\cap b(s\B),\hbox{\ misses\ }
          s\overline\B .\cr}
          \ \ \ \right\}
          \eqno (4.3)
          $$
          
          There is a $\tau_0>0$ such that $\omega \tau_0^2<\delta $ where $\delta$ satisfies (4.3) and
          which is so small that 
          for every $\tau, \ 0<\tau<\tau_0, $ and for every $z\in\R^{M-1}$ the union of all simplices in $\tau(\cS+z)$ 
          contained in $U_0$, contains $U_1$. Let $0<\tau<\tau_0$ and let $z\in\R^{M-1}$. Then $\Gamma (r,\tau, z)\cap W_1$ is a 
          graph over $U_1$ which is contained in $W_1((r-\omega\tau^2,r])$.
          
          We now show that each hyperplane $H$ meeting $\Gamma(r,\tau, z)\cap W$ and tangent to $\Gamma
          (r,\tau, z)$ misses 
          $(r-\omega\tau^2)\overline\B$. This is easy to see. By the convexity of $\Gamma(r,\tau, z)$, all of 
          $\Gamma(r,\tau, z)$ except the simplex $H\cap \Gamma(r,\tau, z)$ is contained in the open halfspace bounded by 
          $H$ which contains the origin. If $H$ would meet $(r-\omega\tau^2)\overline\B$, then, 
          since $\omega\tau^2<\delta$, 
          by (4.3), $H$ would meet $b((r^2-\omega\tau^2)\B)\cap W_1$ at 
          a point not 
          contained in the simplex $H\cap
          \Gamma (r,\tau, z)$, a contradiction.
          
          Dentote by $\pi$ the projection
          $$
          \pi(x_1,\cdots, x_M) = (x_1,\cdots, x_{M-1}).
          $$
          Out of simplices building $\Gamma (r,\tau, z)$ choose those which meet $W$, denote them with 
          $T_j=\Phi_r(S_j), \ 1\leq j\leq n$ and let $C$ be their union. Since $\tau<\tau_0$ 
          the simplices in $\tau(\cS+z)$ contained in 
          $U_0$ cover $U_1$ so the simplices $S_j, 1\leq j\leq n$, cover $U$ and consequently
          $C\cap W$ is a graph over $U$ which cuts $W$ into two connected components. Any path in 
          $W$ connecting points in different 
          components will have to intersect $C$. What remains of $C$ after
          we remove the relative interiors of all simplices $T_j,\ 1\leq j\leq n$, we call 
          the \it skeleton \rm of $C$ and denote by $\hbox{Skel}(C)$. Obviously $\pi (\hbox{Skel}(C))\subset \hbox{Skel}(
          \tau(\cS+z)).$
          
          For each $j, \ 1\leq j\leq n$, there is a linear functional $\ell _j$ on $\C^N$ such that the hyperplane 
          $H_j=\{ z\in\C^N\colon\ \Re(\ell_j(z))=1\} $ contains $T_j$. We know that each $H_j$ misses $(r-\omega\tau^2)
          \overline\B$.
          
          Since $\Gamma(r,\tau, z)$ is convex it is easy to see that given a sufficiently small $\nu>0$ the sets 
          $$
          \{z\in\B\colon \ 1-\nu<\Re \ell_i(z)<1\} ,\ 1\leq i\leq n,
          $$
          intersect $C$ in a small neighbourhood $\cV \subset C$ of $\hbox{Skel}(C)$, and the sets 
          $\{ z\in\B\colon \ \Re \ell_i(z)<1-\nu \} $ contain $(r-\omega\tau^2)\overline\B$. Note that $\cV$ 
          will be arbitrarily small neighbourhood of $\hbox{Skel}(C)$ provided that $\nu >0$ is small enough.
          
          Choose $\varepsilon ,\ 0<\varepsilon<1$. Given $L<\infty $ we use a one variable Runge theorem to 
          get a polynomial $\varphi $ of one variable such that 
          $$ 
          \Re \varphi > L+1\hbox{\ on\ } \{\z\in 2\overline\D\colon \ \Re\z = 1\}, \ \ 
          |\varphi |<\varepsilon /n \hbox{\ on\ } \{\z\in 2\overline\D\colon \ \Re\z <1-\nu\}.
          $$
          Then $\Phi=\sum_{i=1}^n\varphi\circ \ell_i$ is a polynomial on $\C^N$ such that 
          $$
          \Re\Phi >L \hbox{\ on\ }C\setminus \cV,\ \ \ |\Phi|<\varepsilon \hbox{\ on\ } (r-\omega\tau^2)\overline\B .
          $$
          The convex surface $C\subset \hbox{Sh}((r-\omega\tau^2,r])$,\ $C\subset \Gamma (r,\tau, z)$ is such that
          $W\setminus C$ has two components. Let $\cV(\eta )\subset C$ be the $\eta$-neighbourhood of
          $\hbox{Skel}(C)$ where $\eta >0$ is very small. Write $G=W\cap (C\setminus \cV(\eta))$. Then
          $$
          \left.\eqalign{
          (i)\ \ \ &G\subset W((r-\omega\tau^2,r]), \cr
          (ii)\ \  &\hbox{a path\ }p\colon [0,1]\rightarrow W([r-\omega\tau^2,r])
          \setminus G\hbox{\ such that \ }\cr
          & |p(0)|=r-\omega \tau^2,\ |p(1)|=r,\hbox{\ necessarily meets \ }\cV (\eta)\cr
          (iii) \ &\hbox{given\ } L<\infty \hbox{\ and\ } \varepsilon >0\hbox{\ there is a polynomial\ }
          g\hbox{\ on\ }\C^N \hbox{\ such that\ }\cr
          &|g|<\varepsilon \hbox{\ on\ }(r-\omega\tau^2)\overline\B
          \hbox{\ and\ }\Re g>L \hbox{\ on\ } G.\cr}\ \ \ \ \ \right\} 
          \eqno (4.4)
          $$
          Observe also that $\pi$ maps $\cV(\eta)$ to the $\eta$-neighbourhood of $\hbox{Skel}(\tau(\cS + z))$ in $
          \R^{M-1}$.
          \vskip 4mm
          \noindent\bf 5.\ Proof of Lemma 3.1, Part 2. Completion of the proof of Lemma 2.1. \rm
          \vskip 2mm
          Using an easy transversality (or "putting into general position") argument we see that the fact that $\cS$ is periodic with respect to
          $\Lambda$ implies that there are $q_1,\cdots ,\ q_{M-1}\in\R^{M-1}$ such that
          $$
          \hbox{Skel}(\cS)\cap\hbox{Skel}(\cS+q_1)\cap\cdots\cap \hbox{Skel}(\cS+q_{M-1}) =\emptyset 
          $$ which implies that there is a $\mu>0$ such that
          $$
          |x_1-x_0|+|x_2-x_1|+\cdots +|x_{M-1}-x_{M-2}|\geq \mu
          $$ 
          whenever $x_i\in \hbox{Skel}(\cS+q_i)\ \ (0\leq i\leq M-1)$ where  $q_0=0$. 
          It then follows that for every 
          $\tau >0$
          $$
          |x_1-x_0|+|x_2-x_1|+\cdots +|x_{M-1}-x_{M-2}|\geq \tau\mu
          \eqno (5.1)
          $$ 
          whenever $x_i\in\hbox{Skel}(\tau(\cS+q_i)),\ 0\leq i\leq M-1$.
          
          Suppose that $1/2<r_0<r_M<1$. Divide the interval $[r_0,r_M]$ into $M$ equal pieces of length 
          $(r_M-r_0)/M$, let $r_0<r_1<\cdots <r_{M-1}<r_M$ be their endpoints. Choose $\tau,\ 0<\tau<\tau_0$ so small that
          $$
          \tau^2\omega < {{r_M-r_0}\over M}.
          \eqno (5.2)
          $$
          Fix a very small $\eta>0$.  For each $j,\ 1\leq j\leq M$, use
          $\Gamma (r_j,\tau, q_{j-1})$ to construct $G_j\subset \Gamma (r_j, \tau, q_{j-1})$ as in 
          the preceding section. Then for each $j, \ 1\leq j\leq M$, we have $G_j\subset W((r_{j-1},r_j])$ and 
          each path $p\colon [0,1]\rightarrow W([r_{j-1},r_j]),\ |p(0)|=r_{j-1},\  |p(1)|=r_j$ meets the
          $\eta$-neighbourhood of 
          $\hbox{Skel}(\Gamma (r_j,\tau ,q_{j-1}).$ Thus, if $G=\cup_{j=1}^M G_j$ then for each 
          path 
          $p\colon [0,1]\rightarrow W([r_0, r_M])\setminus G$, such that $|p(0)|=r_0,\ |p(1)|=r_M$,  
          there are $t_j,\ 1\leq j\leq M,\ 0<t_1<\cdots <t_M<1$ such that $p(t_j)$ is in the 
          $\eta$-neighbourhood of $\hbox{Skel}(\Gamma
          (r_j, \tau, q_{j-1}))$, so there is a point $z_j\in \hbox{Skel}
          (\Gamma(r_j,\tau, q_{j-1} )$ such that 
          $$
          |z_j-p(t_j)|<\eta\ \ (1\leq j\leq M).
          $$ 
          By (5.1) it follows that
          $$
          \hbox{\ length}(p)\geq \sum_{j=1}^{M-1}|p(t_{j+1})-p(t_j)|\geq \sum_{j=1}^{M-1}|z_{j+1}-z_j| - 
          (M-1)\eta\geq \tau\mu-M\eta.
          $$
          Using (iii) in (4.4) in an induction process again one concludes that given $\varepsilon>0$ and $L<\infty$ 
          there is a polynomial $g$ on $\C^N$ such that $\Re g >L$ on $G$ and $|g|<\varepsilon$ on $r_0\overline\B$.
          
          Thus, if $1/2<r_0<r_M<1$ and if $\tau, \ 0<\tau<\tau_0$ satisfies (5.2) then there is a set 
          $E\subset W((r_0,r_M])$ such that
          
          (i) if $p\colon [0,1]\rightarrow W([r_0,r_M])\setminus E$ is a path such that
          $|p(0)|=r_0,\ |p(1)|=r_M$ then the length of 
          $p$ exceeds $\tau\mu-M\eta$
          
          (ii) given $L<\infty $ and $\varepsilon >0$ there is a polynomial $g$ such that $\Re g >L$ on $E$ 
          and $|g|<\varepsilon$ on $r_0\overline\B$.
          
          We now prove that for every $A<\infty$ and for every segment $J=(\alpha,\beta],\ 1/2<\alpha<\beta<1$ there 
          is a set $E\subset W(J)$ which satisfies (4.1).
          
          So let $1/2<\alpha<\beta<1$. Write 
          $(\alpha,\beta]= J_1\cup J_2\cup\cdots\cup J_\ell$ where
          $J_1<J_2<\cdots <J_\ell$ 
          and where $|J_j|=(\beta-\alpha)/\ell\ \ (1\leq j\leq \ell)$. For each $j$ we shall construct a set 
          $E_j\subset W(J_j)$ as above. 
          Provided that $0<\tau<\tau_0$ and $\tau^2\omega < (\beta - \alpha)/(M\ell)$ the set 
          $E=\cup_{j=1}^\ell E_j$ will then 
          have the property that if a path $p: [0,1]\rightarrow W([\alpha, \beta])\setminus E$\  
          satisfies $|p(0)|=\alpha,\ |p(1)|=\beta $ then
          $$\hbox{length}(p)\geq \ell (\tau\mu-M\eta)= \ell\tau\mu - M\ell\eta .
          $$
          Suppose that $A<\infty $ is given. We show that it is possible to choose $\ell $ and $\tau,\ 0<\tau<\tau_0$, so that
          $$\tau^2\omega < {{\beta-\alpha}\over{M\ell}}
          \eqno (5.3)
          $$
          and 
        $$\ell\tau\mu = A+1 .
        $$
        In fact,  $\tau=(A+1)/(\ell\mu)$ implies that there is an $\ell_0$ such that $0<\tau<\tau_0$ 
        for every $\ell>\ell_0$. For (5.3) to hold we must have
        $$
        \biggl({{A+1}\over{\ell\mu}}\biggr)^2\omega < {{\beta-\alpha}\over {M\ell}}
        $$
        which is clearly possible provided that $\ell>\ell_0$ is chosen large enough. 
        
        So fix such $\ell$ and such $\tau$ and let $\eta >0$ be so small that $M\ell\eta <1$. 
        Then 
        $$
        \hbox{length}(p)\geq \ell\tau\mu-M\ell\eta \geq A+1-1 =A   .
        $$
        Given $L<\infty $ and $\varepsilon >0$, an inductive process again produces a polynomial 
        $\Phi$ such that $\Re \Phi>L$ on $E$ and 
        $|\Phi | <\varepsilon$ on $\alpha\overline\B$. This will complete the proof of 
        Lemma 3.1 and thus the proof of Lemma 2.1 once we have proved that the 
        tessellation $\cS$ can be chosen in such a way that, provided that the ball $U_0$ is small enough, 
        the surfaces $\Gamma (r, \tau, z)$ are convex.
      \vskip 4mm
      \noindent\bf 6.\ Convexity of the surfaces $\Gamma (r,\tau, z)$ \rm 
      \vskip 2mm
      We shall now show how to choose a tessellation $\cS$ in Section 4 so 
      that after choosing $U_0$ small enough the 
      polyhedral surfaces $\Gamma (r,\tau, z)$ will be convex. This is the fact that 
      we used in the proof of 
      Lemma 3.1. We essentially follow [G].
          	  	        
      Perturb the canonical orthonormal basis in $\R^{M-1}$ a little to get an 
       $(M-1)$-tuple of vectors $e_1, e_2, \cdots , e_{M-1}$ in general position so that the lattice
       $$
       \Lambda = \bigl\{ \sum_{i=1}^{M-1}n_ie_i \ \colon \ n_i\in Z, \ 1\leq i\leq M-1\bigr\} 
          	\eqno (6.1)
       $$
        will be generic, and, in particular, no more than $M$ points of $\Lambda$ 
        will lie on the same sphere. 
          		        
        For each point $x\in\Lambda$ there is the\  
        \it Voronei cell \rm\  $V(x)$ consisting of those points of
        $\R^{M-1}$ that are at least as close to $x$ as to any other $y\in\Lambda$, so
        $$
        V(x)=\{ y\in\R^{M-1}\colon\ \hbox{dist} (y,x)\leq 
        \hbox{dist} (y,z)\hbox{\ for all\ }z\in\Lambda
        \} .
        $$
         It is known that the Voronei cells form a tessellation of 
        $R^{M-1}$ and in our case they are all congruent, 
        of the form $V(0)+x, \ x\in\Lambda$ [CS]. 
          		        
        There is a \it Delaunay cell\rm\  for each point that is a vertex of a Voronei cell.
        It is the convex polytope that is the 
        convex hull of the points in $\Lambda $ closest to that point - these points are all on 
        a sphere centered at this point. In our case, when there are no more than $M$ 
        points of $\Lambda$ on a sphere, Delaunay cells are $(M-1)$-simplices. 
        Delaunay cells form a tessellation of $R^{M-1}$ [CS]. In our case it is a \it true \rm 
        Delaunay tessellation , that is, for each cell, the circumsphere of each cell $S$ contains 
        no other points of $\Lambda$ than the vertices of $S$. We shall denote 
        by $\cS$ 
        the family of all simplices - cells of the Delaunay tessellation for the
        lattice $\Lambda $ and this is to be taken as our $\cS$ in Section 4. Clearly $\Lambda$ is
        precisely the set of vertices of the simplices in $\cS$ .
          		        
        The construction implies that the tessellation $\cS$ is periodic with 
        respect to $\Lambda$.
        Thus, there are finitely many simplices $S_1,\cdots, S_\ell$ in $\cS$  such that
         every other simplex of $\cS$ is of the form $S_i + w$ 
         where $w\in\Lambda$ and
          		        $1\leq i\leq \ell$. It is then clear by the periodicity that there 
          		        is an $\eta>0$ such that 
          		        for every simplex $S\in\cS$  
          		        in the $\eta$-neighbourhood 
          		        of the closed ball bounded by the circumsphere 
          		        of $S$ there are no other points of  $\Lambda$ 
          		        than the vertices of $S$. 
          		        
          		        When we pass from $\cS$ to $\tau(\cS+z)$ where $\tau>0$ and $z\in\R^{M-1}$ everything  
          		        changes proportionally. For instance, for every simplex $S\in \tau(\cS+z)$ in the 
          		        $(\tau\eta)$-neighbourhood of the closed ball bounded by the circumsphere of $S$ there will be 
          		        no other vertex of $\tau(\cS +z)$ than the vertices of $S$.
          		        
          		        We must now show that if $U_0$ is chosen sufficiently small then for every $r,\ 1/2<r<1$, 
          		        every $\tau>0$ and every $z\in\R^{M-1}$ for every simplex $S\in\tau(\cS+z)$ contained in 
          		        $U_0$ the intersection of the hyperplane $H$ containing $\Psi_r(S)$ with $\Gamma(r,\tau, z)$ is
          		        precisely $\Psi_r(S)$. 
          		        
          		        So let $S$ be such a simplex and let $H$ be the hyperplane in $\R^M$ containing $\Psi_r(S)$. Then
          		         $H$ intersects $b(r\B)$ in a sphere $\Gamma$ that is the circumsphere of $\Psi_r(S)$ in $H$. 
          		         One component of $b(r\B)\setminus\Gamma$ is contained in the open halfspace bounded by $H$ 
          		         which contains the origin. 
          		         Denote this component of $b(r\B)\setminus\Gamma$ by $E$.
          		         
                                Recall that all vertices of $\Gamma(r,\tau,z)$ are contained in $b(r\B)$. Thus, to prove that 
                                $H$ contains no other vertex of $\Gamma(r,\tau, z)$ than the vertices of $\Psi_r(S)$ it is enough to show that
         $$
         \hbox{all vertices of\ }\Gamma(r,\tau,z)\hbox{\ except the vertices of \ }\Psi_r(S)
         \hbox{\ are contained in\ }E.
         \eqno (6.2)
         $$
         Since $\pi|W_0\cap b(r\B)\colon\ W_0\cap b(r\B )\rightarrow U_0$ is one to one, 
         to show (6.2) it is enough to show that
         $$
         \left.\eqalign{
         &\hbox{the vertices of all simplices\ }T\in \tau(\cS+z)\
         \hbox{contained in\ } U_0 \cr
         &\hbox {except the vertices of\ }
         S\hbox{\ are  contained}\hbox{ in the complement}\cr
         &\hbox{of the bounded domain in\ } \R^{M-1}
         \hbox{\ bounded by \ }\pi(\Gamma ).\cr}\ \ 
         \right\}
         \eqno (6.3)
         $$
         Since the $(\tau\eta)$-neighbourhood of the closed ball in $\R^{M-1}$ bounded by 
         the circumsphere of $S$ contains no other vertex 
         of $\tau(\cS+z)$ than the vertices of $S$ it follows that to show (6.3) it is enough 
         to show that
         $$
         \left. 
         \eqalign{
         &\hbox{provided that \ } U_0  \hbox{\ is chosen small enough 
         on the outset then for every \ }
         \cr
         &S\in\tau(\cS+z)\hbox{\ contained\ }
         \hbox{\ in\ }U_0\hbox{\ the projection\ }  \pi(\Gamma) 
         \hbox{\ of the circumsphere\ }\cr
         &\Gamma\hbox{\ of\ }\Psi_r(S)
         \hbox{\ in the hyperplane\ }H\hbox{\ containing\ }
         \Psi_r(S)\hbox{\ is contained in the }\cr
         &(\tau\eta)\hbox{-neighbourhood
         of the circumsphere of\ } S 
         \hbox{\ in \ }\R^{M-1} .\cr}
         \right \}
         \eqno (6.4)
         $$
Given a $(M-1)$-simplex $T\subset \R^M$ denote by $\Gamma (T)$ the circumsphere of $T$ in 
the hyperplane containing $T$. Given a $(M-1)$-simplex $S\subset \R^{M-1}$ with vertices $v_1,
\cdots ,v_M $ and $\omega>0$ denote by $\Omega_\omega (S)$ the set of all simplices with 
vertices $(v_1,q_1),\cdots , (v_M,q_M)$ where $q_i\in\R $ satisfy 
$$
|q_i-q_j|\leq \omega |v_i-v_j|\ \ (1\leq i, j\leq M) .
$$
\vskip 2mm \noindent \bf PROPOSITION 6.1\ \ \it Let $S\subset \R^{M-1}$ be a $(M-1)$-simplex. 
Given $\eta >0$ there is an $\omega >0$ such that for every $T\in\Omega_\omega(S)$ the 
projection $\pi(\Gamma (T))$ is contained in the $\eta$-neighbourhood of $\Gamma (S)$. Moreover, 
for any $\tau >0$ and for any $T\in\Omega_\omega (\tau S)$ the projection
$\pi (\Gamma (T))$ is contained in 
the $(\tau\eta)$-neighbourhood of $\Gamma (\tau S)$. \rm
\vskip 1mm\noindent\bf Proof.\ \ \rm Let $S\subset \R^{M-1}$ be a simplex with 
vertices $v_1,\cdots, v_M$ and let $T$ be a simplex with vertices
$(v_1,q_1),\cdots (v_M,q_M)$. Note that $\pi(\Gamma (T))$ does not change if we 
translate $T$ in the 
direction of the last axis so with no loss of generality consider  the simplex with the 
vertices $(v_1,q_1-q_M),\cdots (v_{M-1},q_{M-1}-q_M), (v_M,0)$. We now show that 
if $P$ is a simplex with 
vertices $w_1=(v_1,\beta_1),\cdots ,w_{M-1}= (v_{M-1},\beta_{M-1}),  w_M=(v_M,0)$ then
$$\left.\eqalign{
\pi(\Gamma(P))\hbox{\ is arbitrarily close to \ }\Gamma(S)\hbox{\ provided that\ }
\beta_1,\cdots ,\beta_{M-1}& \cr 
\hbox{\ are sufficiently small.\ }& \cr}\ \right\}
\eqno (6.5)
$$
This implies that if $\eta >0$ then there is an $\varepsilon >0$ such that if $|q_i-q_M|
<\varepsilon \ (1\leq i\leq M-1)$ then $\pi (\Gamma (T))$
is contained in the $\eta$-neighbourhood of $\Gamma (S)$. Picking now $\omega >0$ 
so small that $\omega |v_i-v_M|<\varepsilon\ \ (1\leq i\leq M-1)$
 completes the proof of the first part of proposition. To prove (6.5), let $H$ be the hyperplane 
containing $P$ and for each $j,\ 1\leq j\leq M-1$, let $H_j$ be 
the hyperplane through the midpoint of the segment with endpoints $w_j,w_M$, perpendicular to 
$w_M-w_j$. The center $C$ of $\Gamma (P)$ is the 
intersection of $H, H_1,\cdots, H_{M-1}$. Since these hyperplanes are 
in general position and change continuously with 
$\beta_1,\cdots ,\beta_{M-1}$, the point $C$ and consequently $\Gamma (P)$ changes 
continuously with $\beta_1,\cdots,\beta _{M-1}$. When $\beta_1=\cdots =\beta_{M-1}=0$ we 
have $P=S$ so $\Gamma (P)=\pi (\Gamma(P))=\Gamma (S)$. This implies (6.5).

To prove the second statement of the proposition assume that $\tau>0$ and that 
$T\in\Omega_\omega (\tau S)$. Then the vertices of $T$ are $(\tau v_1,p_1), \cdots (\tau v_M, p_M)
$ where $|p_i-p_j|\leq \omega |\tau v_i-\tau v_j|\ \ (1\leq i, j\leq M)$. Writing $p_i=\tau q_i$ we get that 
$$
|q_i-q_j|\leq\omega |v_i-v_j|\ \ (1\leq i, j\leq M)
\eqno (6.6)
$$
Thus, the vertices of $T$ are $(\tau v_1,\tau q_1),\cdots (\tau v_M,\tau q_M)$ where  (6.6) holds, 
that is $T=\tau \tilde S$ where $\tilde S\in\Omega_\omega (S)$. 
Clearly $\Gamma (T)=\tau\Gamma (\tilde S)$ and so $\pi (\Gamma (T)) = \pi 
(\tau\Gamma (\tilde S)) = \tau \pi (\Gamma (\tilde S)).$ Since $\tilde S\in\Omega_\omega (S)$ the 
preceding discussion shows that 
$\pi (\Gamma (\tilde S))$ is contained in the $\eta$-neighbourhood of $\pi (\Gamma (S))$ so it follows that $\pi (\Gamma (T))$ 
is contained in the $(\tau\eta)$neighbourhood of $\Gamma(\tau S)$. This completes the 
proof of Proposition 6.1.

To prove (6.4) recall first that there are finitely many simplices $S_1,\cdots, S_\ell$ in $\cS$  such
that every other simplex of $\cS$ is of the form $S_i + w$ where $w\in\Lambda$ and 
$1\leq i\leq \ell$. Thus, there is an $\omega >0$ such that the statement of Proposition 6.1 
holds for every simplex 
$S\in \cS+z$. Recall that $\hbox{grad}\psi_r$ vanishes at the origin so one can choose $U_0$, a ball 
centered at the origin, so small that
$$
|(\hbox{grad}\psi_r)(x)|<\omega\hbox{\ for all\ }x\in U_0\hbox{\ and all\ }r, \ 1/2<r<1 .
\eqno (6.7)
$$
This implies that for every $S\in\tau(\cS + z)$ with vertices $v_1, v_2,\cdots, v_M$, contained in $U_0$, 
the simplex $\Psi_r(S)$ with vertices $(v_1,\psi_r(v_1)),\cdots ,(v_M,\psi_r(v_M))$, by (6.7), satisfies 
$$
|\psi_r(v_i)-\psi_r(v_j)| \leq \omega |v_i-v_j|\ \ (1\leq i, j\leq M)
$$
so the simplex $\Psi_r(S)$ 
belongs to $\Omega_\omega (S)$ so
(6.4) follows by Proposition 6.1. This completes the proof of convexity of surfaces 
$\Gamma(r,\tau, z)$ and completes the proof of Lemma 3.1. 
The proof of Lemma 2.1 is complete. Theorem 1.1 is proved.

          		        \vskip 10mm
          		        This work was supported by the Research Program P1-0291 from ARRS, Republic of Slovenia.
				        \vfill
				        \eject
				        \centerline{\bf REFERENCES}
				        \vskip 4mm
				        \rm
				        
				 \noindent [AL]\ A.\ Alarc\'{o}n and F.\ J.\ L\'{o}pez:\ Complete bounded complex curves in $\C^2$.
				        
				        \noindent To appear in J.\ Europ.\ Math.\ Soc.  http://arxiv.org/pdf/1305.2118.pdf
				        \vskip 1mm
				        \noindent [B]\ A.\ Brondsted:\ \it An Introduction to Convex Polytopes.\rm
				        
				        \noindent Graduate Texts in Math.\ 90, 1983.  Springer-Verlag New York Inc.
				        \vskip 1mm
				        \noindent [CS]\ J.\ H.\ Conway, N.\ J.\ A.\ Sloane:\ \it Sphere Packings, Lattices 
				        and Groups. \rm 
				        
				        \noindent Grundl.\ Math.\ Wiss. 290, 1988 Springer Verlag New York Inc.
				        \vskip 1mm
				        \noindent [G]\ J.\ Globevnik:\ A complete complex surface in the ball of $\C^N$.
				        
				        \noindent Submitted for publication, January 2014. http://arxiv.org/pdf/1401.3135.pdf
				        \vskip 1mm
				        \noindent [GS]\ J.\ Globevnik and E.\ L.\ Stout:\ Holomorphic functions with highly 
				        noncontinuable boundary behavior.
				        
				        \noindent J.\ Anal.\ Math.\ 41 (1982) 211-216
				        \vskip 1mm
				        \noindent [H]\ L.\ Hormander:\ \it An Introduction to Complex Analysis in Several
				        variables. \rm
				        
				        \noindent North Holland, Amsterdam-London, 1973
				        \vskip 1mm
				        \noindent [J]\ P.\ W.\ Jones:\ A complete bounded complex submanifold of $\C^3$.
				        
				        \noindent Proc.\ Amer.\ Math.\ Soc.\ 76 (1979) 305-306
				        \vskip 1mm
				        \noindent [Y1]\ P.\ Yang:\ Curvature of complex submanifolds of $\C^n$.
				        
				        \noindent J.\ Diff.\ Geom.\ 12 (1977) 499-511
				        \vskip 1mm
				        \noindent [Y2]\ P.\ Yang:\ Curvature of complex submanifolds of $\C^n$. 
				        
				        \noindent In: Proc.\ Symp.\ Pure.\ Math.\ Vol.\ 30, part 2, pp. 135-137. Amer.\ Math.\ Soc., 
				        Providence, R.\ I.\ 1977
				        \vskip 10mm
				        \noindent Institute of Mathematics, Physics and Mechanics
				        
				        \noindent Ljubljana, Slovenia
				        
				        \noindent josip.globevnik@fmf.uni-lj.si
				        
				        \end